%% file: _2.tex
\documentclass[winedt,
                       yap
                            ]{pit}
\usepackage{srcltx}
\usepackage{amsmath}
\usepackage{amssymb}
\usepackage{array}
\usepackage{longtable}
\usepackage{tabularx}
\usepackage{epsfig}
\usepackage{psfrag}
\usepackage{enumerate}
\usepackage{url}
\usepackage{mathrsfs}
\usepackage{dsfont} % instead of bbm !!!

\usepackage[mathscr]{euscript}

\makeatletter
\newcommand{\mylabel}[2]{#2\def\@currentlabel{#2}\label{#1}}
\makeatother

\usepackage{tikz}
\usetikzlibrary{arrows}
\tikzstyle{int}=[draw, fill=white!20, minimum size=2em] \tikzstyle{init}=[pin
edge={to-,thin,black}]

\input{newsymb.tex}

\allowdisplaybreaks

\begin{document}
\frontmatter % for the preliminaries
%
%Here may be written anything what you want to be included BEFORE
%the table of contents and the main body of the journal issue
%
%
%
%Below variables should be specified for EACH issue of journal
%
\JournalName{Problems of Information Transmission}% by default -- Problems of Information Transmission
\JournalISSNCode{0032-9460}
\TransYearOfIssue{2021}
\TransCopyrightYear{2021}
\TransVolumeNo{57}
\TransIssueNo{2}
\OrigYearOfIssue{2021}
\OrigCopyrightYear{2021}
\OrigVolumeNo{57}
\OrigIssueNo{2}

%\renewcommand{\baselinestretch}{1.05}
%\tableofcontents % <--- use when preparing the final version
%
%\renewcommand{\baselinestretch}{1}
\mainmatter % start of the contributions
\setcounter{page}{99}
%
%Below variables should be specified for EACH ARTICLE before the %\input command
%or inside the corresponding input article file
%
%\title{}% OBLIGATORY!
%\author{}% OBLIGATORY!
%\institute{}% OBLIGATORY!
%\received{}% OBLIGATORY!
%\titlerunning{}
%\authorrunning{}
%\toctitle{}
%\tocauthor{}
%\OrigPages{}% OBLIGATORY!
%\OrigCopyrightedAuthors{}% OBLIGATORY!

\Rubrika{Information Theory}% The Rubrika name if necessary

\Rubrika{Coding Theory}

\Rubrika{Large Systems}

\input{log-mog.tex}

\Rubrika{Information Protection}

%\input{kaba.tex}

%\Rubrika{Automata Theory}

%\Rubrika{Source Coding}

%\Rubrika{Retraction Notes}

%\Rubrika{Methods of Signal Processing}

%\Rubrika{Communication Network Theory}

%\Rubrika{The International Dobrushin Prize}

%\Rubrika{Information Theory and Coding Theory}

%\Rubrika{Image Recognition}

%\Rubrika{Letters to the Editor}

%\Rubrika{Chronicles}

%\Rubrika{Obituary}

%\Rubrika{}
%
%\input{auth-ind.tex}
%
%\input{contents.tex}

\end{document}

%% file: newsymb.tex
\def\bfm#1{\boldsymbol{#1}}

\let\kappa\varkappa

\let\epsilon\varepsilon

\let\phi\varphi
\let\upn=\textup

\def\etal{{\noexpand\lowercase{\noexpand\normalfont et al.}}}
\newcolumntype{C}{>{\centering\arraybackslash}X}
\def\E{\operatorname{\mathbf E}\nolimits}
\def\P{\operatorname{\mathbf P}\nolimits}

\def\nfrac#1#2{\mbox{\small$\dfrac{#1}{#2}$}}

\def\mmid{\mathchoice{\hskip1.5pt|\hskip1.5pt}{\hskip1.5pt|\hskip1.5pt}{\hskip0.5pt|\hskip0.5pt}{\hskip0.3pt|\hskip0.3pt}}

\makeatletter
\renewcommand{\pod}[1]{\if@display\mkern10mu\else\mkern6mu\fi(#1)}
\renewcommand{\pmod}[1]{\pod{{\operator@font mod}\mkern6mu#1}}
\makeatother
\renewcommand{\int}{\intop\limits}
\def\?%
{}

%% file: log-mog.tex
\def\essinf{\operatorname{ess\,inf}\limits}
\def\esssup{\operatorname{ess\,sup}\limits}

\OrigPages{71--89}

\OrigCopyrightedAuthors{The Author(s)}

\title{Limit Theorems for the Maximum Path Weight \\ in a Directed Graph on the Line with
Random Weights of Edges
%\footnote{The paper is supported by the Mathematical Center in Akademgorodok, grant 075-15-2019-1675 by the Ministry of Science and Higher Education; and also by partial support from RFBR grant 18-11-00129 (A.V. Logachov and A.A. Mogulskii) and joint Russian-French grant RFBR 19-51-15001 -- CNRS 193-382 (S.G. Foss and T. Константопулос).}
}

\titlerunning{Limit Theorems for the Maximum Weight}
\authorrunning{Konstantopoulos \etal}

\received{Received November 19, 2020; revised February 1, 2021; accepted February 8,
2021}

\author{T. Konstantopoulos\inst{1}\zpt\dovesok{$^*$}, A. V.
Logachov\inst{2}\zpt\inst{3}\zpt\inst{4}\zpt\dovesok{$^{**}$}, A. A.
Mogulskii\inst{2}\zpt\inst{3}\zpt\dovesok{$^{***}$},\\ and S. G.
Foss\inst{2}\zpt\inst{3}\zpt\inst{5}\zpt\dovesok{$^{****}$}}

\institute{Department of Mathematical Sciences, University of Liverpool, Liverpool,
UK \and Sobolev Institute of Mathematics, Siberian Branch\\ of the Russian Academy of
Sciences, Novosibirsk, Russia \and Novosibirsk State University, Novosibirsk, Russia
\and Siberian State University of Geosystems and Technologies \and School of
Mathematical Sciences, Heriot--Watt University, Edinburgh, UK\\
\email{$^*$T.Konstantopoulos@liverpool.ac.uk}\textup,
\rm$^{**}$\texttt{omboldovskaya@mail.ru}\textup,\\
\rm$^{***}$\texttt{mogul@math.nsc.ru}\textup,
\rm$^{****}$\texttt{sergueiorfoss25@gmail.com}}

\doi{210200??}

\renewcommand{\aboveinst}{12pt}
\maketitle
\renewcommand{\aboveinst}{6pt}

\renewcommand\theequation{\arabic{section}.\arabic{equation}}
\numberwithin{equation}{section}

\begin{abstract}
We consider an infinite directed graph with vertices numbered by integers $ \ldots, -2,-1,0,1,2, \ldots$, where any pair of vertices $j<k$ is connected by an edge $(j,k)$ that is directed from $j$ to $k$ and has a random weight $v_{j,k}\in [-\infty,\infty)$. Here $\{v_{j, k}, j<k\}$ is a family of independent and identically distributed random variables that take either finite values (of any sign) or the value $ - \infty$.
A path in the graph is a sequence of connected edges $(j_0, j_1), (j_1,j_2),\ldots, (j_{m-1}, j_m)$ (where $j_0<j_1<\ldots <j_m$), and its weight is the sum 
$\sum_{s=1}^{m} v_{j_{s-1},j_s}\ge -\infty$ 
of the weights of the edges. Let $w_{0,n}$ be the maximal weight of all paths from $0$ to $n$. \\
Assuming that ${\mathbf P}(v_{0,1}>0)>0$, that the conditional distribution ${\mathbf P} (v_{0,1}\in \cdot \ | \ v_{0,1}>0 )$ is nondegenerate, and that ${\mathbf E} \exp (Cv_{0,1})<\infty$ for
some $C=const >0$, we
study the asymptotic behaviour of the random sequence $w_{0, n}$ as $n\to\infty$.
In the domain of the normal and moderately large deviations we obtain a local limit theorem when the distribution of random variables $v_{i,j}$ is arithmetic and an integro-local limit theorem if this distribution is non-lattice.
\end{abstract}

\section{Introduction, main notation and the main result}

We consider an infinite directed graph
$G(\mathbb{Z},E)$, with vertices indexed by all integers $\mathbb{Z}=\{\ldots,-2,-1,0,1,2,\ldots\}$,
whose edges  $E=\{e=(j,k),\;
j<k,\: j,k\in\mathbb{Z}\}$ are all edges directed from smaller to bigger vertices. We assume that there no directed edges from bigger to smaller vertices and that there is no loops of the form $(j,j)$.

Every edge
$e\in E$ gets a weight 
$v_e$,
that may be either a number (positive or negative) or 
$-\infty$.
We assume that the random variables $\{ v_{j,k}, j<k\}$ are mutually independent and distributed as a random variable $v$ taking values in $[-\infty,\infty)$.
Let $p=\P(v>-\infty)$ and $p^+=\P(v>0)$.
Let $v^+$ be a random variable with distribution 
\begin{align}\label{vplus}
\P(v^+<t)=\P(v< t\mmid v>0), \quad t>0.
\end{align}
Throughout the paper, we assume that the following conditions hold:
\begin{align}\label{maincond}
p^+>0, \qquad \P(v^+=c)<1 \quad \text{for any}\  c>0,\qquad  
\E e^{Cv^+}<\infty\quad \text{for some}\ C>0,
\end{align}
i.e. the random variable $v$ takes positive values with positive probability, its distribution on the positive halfline is non-degenerate, and the right tail of its distribution is relatively light.

In the literature, one can also find another interpretation of the model: if 
$v_e=-\infty$, then one can say that there is no edge
$e$,
and if
$v_e>-\infty$,
then the edge exists and its weight equals 
$v_e$. In the latter setting, we get two independent ``randomnesses'': an edge may either exist or not, and if it exists, then its weight is an independent of everything random variable with distribution $\P(v\in\cdot \mmid v>-\infty)$.

A \emph{path} $\pi$ of length $L(\pi)=m$ is a sequence of $m$ connected edges $e_1=(j_0,j_1)$, $e_2=(j_1,j_2),\ldots,e_m=(j_{m-1},j_m)$ where the end vertex of each edge coincides with the initial vertex of the next edge and $j_0<j_1<\ldots<j_m$, and we say that this is a path from
$j_0$ to $j_m$  and write $e_i\in\pi$,
$i=1,\ldots,L(\pi)$. 
The weight $w(\pi)$ of the path is defined as the sum of the weights of its edges, i.e.
$$
w(\pi)=\sum\limits_{s=1}^{L(\pi)}v_{j_{s-1},j_s}=\sum\limits_{e\in\pi}v_e.
$$
Clearly, the path weight is finite if and only if all weights of its edges are finite.

For $j<k$, let $\Pi_{j,k}$ be the family of all paths from $j$ to $k$ having finite weights
(i.e.
 $w(\pi)>-\infty$ for all $\pi\in \Pi_{j,k}$)
and let  $w_{j,k}$ be the maximal weight of all paths from $j$ to $k$. Then
$$
w_{j,k}=\max\limits_{\pi\in\Pi_{j,k}} w(\pi)
$$ 
with probability 1, since we follow the standard convention that the maximum over empty set is
 $-\infty$.
 We also let $w_{j,j}=0$, for all $j$.

Such graphs with random weights naturally appear in various applications. For example, if the edge weight takes only two values, $1$ (the edge exists) or $-\infty$ (no edge), i.e.
\begin{align}\label{1-infty}
p={\mathbf P} (v=1) =
1 - {\mathbf P} (v=-\infty),
\end{align}
then such a graph may describe ordering of jobs in a computer network (see, e.g. \cite{Coh, New1}),
where vertices represent jobs and edges their time constraints (if $v_{j,k}=1$, then service of job $k$ cannot start before service of job $j$ ends); or functioning of biological models (see, e.g., \cite{Gel, New2}) where vertices represent types of animals and paths describe "food chains": 
if $v_{j,k}=1$, then type $k$ may be considered as a food for type $j$).

Introduce two mutually exclusive conditions:
\begin{itemize}
\item[$\lbrack\mathbf{R}\rbrack$]\vskip-5pt
\emph{The distribution of the random variable $v$ is non-lattice,} this means that, for any $a$ and $h>0$, the probability that $v$ takes values in the lattice of span $h$ shifted by $a$ is strictly smaller than 1, $\sum\limits_{s=-\infty}^{\infty}\P(v=a+sh)<1$.
\item[$\lbrack\mathbf{Z}\rbrack$]
\emph{The distribution of random variable $v$ is arithmetic}, i.e.
$\sum\limits_{s=-\infty}^{\infty} \P(v=sh)=1$, for some $h>0$.
Without loss of generality, we may assume further that the lattice span is $h=1$, this means that $v$  is an integer-valued random variable and the greatest common divisor of the set 
$\{k\ge 1: \P(v=k)>0\}$ equals one.\footnote{
Conditions $[\mathbf{R}]$ and $[\mathbf{Z}]$ may be also formulated in terms of the characteristic function $f(z)=\mathbf{E}e^{izv}$ of random variable $v$. Namely,\\
$[\mathbf{R}]$   $|f(2\pi z)|<1$ for all $z \neq 0$;\\
$[\mathbf{Z}]$
$f(2\pi z)=1$ for all $z\in \mathbb{Z}$ and
$|f(2\pi z)|<1$ for all $z\in \mathbb{R}\setminus\mathbb{Z}$.
}.
\end{itemize}\vskip-5pt

Notice that we have excluded the case $v^+=const$ (see \eqref{maincond}) and the case of lattice, but non-arithmetic distributions.

We are interested in the asymptotic behaviour of the random sequence
 $w_{0,n}$ as $n\rightarrow\infty$.
 We will consider the domain of the normal and moderately large deviations and prove
 a local limit theorem under condition 
 $[\mathbf{Z}]$ and an integro-local limit theorem under condition $[\mathbf{R}]$.

The proof of these statements is split into two steps. In the first step (Section~2) we first introduce an embedded regenerative sequence with corresponding ``weights''
(using the methods developed in the papers \cite{Foss3,Foss2})
and then show that both the lengths of regenerative cycles $\tau_k$ and the corresponding cycle weights $\zeta_k$ have finite exponential moments (precise definitions of these variables are given in Section~3). Note that the sequence $\left\{(\tau_k,\zeta_k)\right\}_{k=1}^{\infty}$ contains  independent, for $k\ge 1$, and identically distributed, for $k\geq 2$, two-dimensional random vectors that have, for $k\ge 2$, a common distribution with a random vector $(\tau,\zeta)$, whose coordinates do typically depend on each other.

In the second step of the proof (Section~3), we note that vectors 
 $(\tau_k,\zeta_k)$ form a stationary compound renewal process (CRP), and one can apply to its study methods and results from the paper 
\cite{Mog1}. Then we prove that, in the limit theorems, the asymptotics of the sequence $w_{0,n}$
coincides with that of the introduced CRP, and this completes the proof of our results. 

In order to formulate our results, it is left to introduce the rate function for the CRP driven by random vector 
$(\tau,\zeta)\stackrel{d}{=} (\tau_2,\zeta_2)$.
For $(\lambda,\mu)\in\mathbb{R}^2$, let 
\begin{equation}\label{17.11.1}
A(\lambda,\mu):=\ln {\bf E}e^{\lambda\tau+\mu\zeta}.
\end{equation}
Introduce the convex set 
$$
\mathcal{A}^{\le 0}:=\{(\lambda,\mu):\: A(\lambda,\mu)\le 0\}
$$
and let 
$$
D(\alpha):=\sup_{(\lambda,\mu)\in \mathcal{A}^{\le 0}}\{\lambda+\mu\alpha\}.
$$
The rate function $D(\alpha)$ plays a certain role in the description of the logarithmic asymptotics of the large deviations probabilities for the CRP determined by vector $(\tau,\zeta)$,
and it has been studied quite thoroughly (see, e.g., \cite{Mog2}). Notice that this is a convex non-negative function that takes value $0$ at only one point 
 $\alpha=a$ where 
\begin{equation}\label{17.11.2}
a= \frac{\E\zeta}{\E\tau}>0.
\end{equation}
Under our assumptions, function $D(\alpha)$ is analytic in a neighbourhood of point $\alpha=a$
and, further, 
$$
D(a)=0,\qquad D'(a)=0,\qquad D''(a)= \frac{1}{\sigma^2},
$$
where
\begin{equation}\label{17.11.3}
\sigma^2:=\frac{{\bf E}(\zeta-a\tau)^2}{\E\tau}.
\end{equation}

We present now the main result of the paper.

\begin{theorem}\label{th2}
Assume that conditions \eqref{maincond} hold.
\begin{enumerate}[\rm I.]
\item\vskip-5pt
If the random variable
$v$ satisfies conditions $[\mathbf{Z}]$, then, for any sequence $x=x_n\in\mathbb{Z}$ such that  $\alpha:={x}/{n}\to a$ as 
$n\to\infty$, the following asymptotic relation holds:
\begin{equation}\label{1.3}
\P(w_{0,n}=x) \sim \frac{1}{\sigma\sqrt{2\pi n}} e^{-nD(\alpha)}.
\end{equation}
If, in addition, 
$y_n:=x-an=o(n^{2/3})$, then we get  
$$
\P(w_{0,n}=x) \sim \frac{1}{\sigma\sqrt{2\pi n}} e^{-\frac{y_n^2}{2n\sigma^2}}.
$$
\item
If the random variable
$v$ satisfies conditions $[\mathbf{R}]$, then, for a certain sequence of positive numbers $\Delta^{(0)}_n=o(1)$ and for any sequence $x=x_n\in\mathbb{R}$ such that
$\alpha:={x}/{n}\to a$ $n\to\infty$, the following asymptotic relation holds: 
\begin{equation}\label{1.2}
\P(w_{0,n}\in [x,x+\Delta_n)) \sim \frac{\Delta_n}{\sigma\sqrt{2\pi n}}
e^{-nD(\alpha)},
\end{equation}
where sequence
$\Delta_n=o(1)$ satisfies relations $\Delta_n\ge\Delta^{(0)}_n$ (i.e. converges to $0$
sufficiently slowly).
If, in addition,
$y_n:=x-an=o(n^{2/3})$, then we get 
$$
\P(w_{0,n}\in [x,x+\Delta_n)) \sim \frac{\Delta_n}{\sigma\sqrt{2\pi n}}
e^{-\frac{y_n^2}{2n\sigma^2}}.
$$
\end{enumerate}
\end{theorem}

\begin{remark}
One can strengthen the results of Theorem  
\ref{th2} by considering, along with the normal and moderately large deviations of the type $\alpha:={x}/{n}\to a$, the large deviations of the type $|\alpha-a|\le\delta$ for some (generally speaking, small) $\delta>0$. In this case the constant $\nfrac{1}{\sigma\sqrt{2\pi}}$ that appears in the right-hand sides of relations \eqref{1.3},
\eqref{1.2} should be replaced by a more complicated function that depends on parameter $\alpha={x}/{n}$. However, for determining this function, one needs to produce additional complicated constructions. This is why we have decided to restrict our consideration in  Theorem \ref{th2} to the normal and moderately large deviations.
\end{remark}

\begin{remark}
The statements of Theorem \ref{th2} are presented in terms of rate function $D(\alpha)$ that is determined by the distribution of random vector $(\tau,\zeta)$ that is introduced implicitly and depends on parameters $c_1$ and $c_2$ that are chosen arbitrarily from a certain interval (see Lemma  \ref{l1}). In Theorem~3 below, we will show that the results of  Theorem \ref{th2} do not depend on a particular choice of these parameters.
\end{remark}

The asymptotic properties of sequence
$w_{0,n}$, as $n\to\infty$, have been studied earlier in the papers \cite{Foss3,Foss2,Foss1,Foss4}. In the paper \cite{Foss1}, the authors consider the case $\P(v>0)=1$ and prove the strong law of large numbers and the central limit theorem assuming that the third moment of random variable $v$ is finite, and limit theorems of another type if the latter condition is violated. 
The central limit theorem is the case of signed random variable $v$ is proved in \cite{Foss4}. In the earlier papers \cite{Foss3,Foss2}, the case of weights \eqref{1-infty} was studied. We should mention also the paper \cite{Tes} where the asymptotics for the \emph{minimal} path length from $0$ to $n$
were considered, as 
$n\to\infty$, in the case where the weights are constant, but the probabilities of existence of edges depend on distances between the vertices.  

The rest of the paper includes three Sections. In Sections~2 and 3 we provide the proofs of our results following the scheme presented above, and Section~4 contains an auxiliary result.

\section{Construction of regenerative sequence and its properties}

In this Section, we introduce a construction that allows us to determine an a.s. infinite random set  of vertices $\left\{\Gamma_i\right\}_{i\in\mathbb{Z}}$ (in what follows, we call them \emph{renewal vertices} -- see Definition~3 below), where 
$\ldots<\Gamma_{-2}<\Gamma_{-1}<0\le\Gamma_0 <\Gamma_1 <\ldots\strut$, such that:
\begin{enumerate}
\item\vskip-5pt
Sequence of two-dimensional vectors
\begin{equation}\label{tauzeta1}
(\Gamma_n-\Gamma_{n-1}, w_{\Gamma_{n-1},\Gamma_n}),\quad n\ne 0,
\end{equation}
consists of independent and identically distributed vectors that do not depend on 
$(\Gamma_{-1},\Gamma_0, w_{\Gamma_{-1},0},w_{0,\Gamma_0})$. Using the terminology of the theory of point processes, one can say that sequence
$\{(\Gamma_n,w_{\Gamma_{n-1},\Gamma_{n}})\}$ forms a stationary marked point process with marks $\{w_{\Gamma_{n-1},\Gamma_{n}}\}$, that determines the CRP.
\item\vskip2pt
For some 
$C>0$, all four exponential moments
\begin{equation}\label{tauzeta2}
\E\exp (C\Gamma_0),\qquad \E\exp (C(\Gamma_1-\Gamma_0)),\qquad \E\exp
(Cw_{0,\Gamma_0}),\qquad \E\exp (Cw_{\Gamma_0,\Gamma_1})
\end{equation}
are finite. Then, for $C_1=C/2$, the following moments are finite with necessity, too:
\begin{equation}\label{tauzeta3}
\E\exp (C_1(\Gamma_0+w_{0,\Gamma_0}))\qquad \text{and}\qquad \E\exp
(C_1(\Gamma_1-\Gamma_0+w_{\Gamma_0,\Gamma_1})).
\end{equation}
\item
For any 
$0\le m\le n$, if $\Gamma_m\le n$, then 
\begin{equation}\label{mainpres}
w_{0,n}=w_{0,\Gamma_0}+w_{\Gamma_0,\Gamma_1}+\ldots+w_{\Gamma_{m-1},\Gamma_m}+
w_{\Gamma_m,n}
\end{equation}
\end{enumerate}\vskip-5pt
(recall that we assume $w_{j,j}=0$, for $j\in\mathbb{Z}$).

We need a number of auxiliary statements. We will partially follow the scheme of the proof of one of the main results in \cite{Foss1} where the weights were assumed to take either positive values or value $-\infty$ and where conditions for existence of the first and second moments of
random variables $\Gamma_0$ and $w_{0,\Gamma_0}$ were studied.

\subsection{Construction of skeleton and renewal points}

We will introduce consequently four random subsets of the set
$\mathbb{Z}$ of vertices: the set $\mathcal{S}$ of skeleton, the set $\mathcal{S}^+$ of skeleton-plus, the set $\mathcal{R}$ of renewal, and the set $\mathcal{R}^+$ of renewal-plus vertices.

\begin{definition} Vertex
$x\in\mathbb{Z}$ is called a \emph{skeleton} vertex if it is connected to any other vertex by a path of finite weight, i.e. for any $j<x$ and $k>x$ inequalities $w_{j,x}>-\infty$ and  $w_{x,k}>-\infty$ hold. Denote by 
$\mathcal{S}$ the random set of skeleton vertices.
\end{definition}

If $p=\P(v>-\infty)=1$, then every vertex $x\in\mathbb{Z}$ is a skeleton vertex. If $p\in (0,1)$, then Lemmas 5--7 from \cite{Foss2} imply the following five statements.
\begin{enumerate}
\item\vskip-5pt
The probability for vertex 
$x$ to be a skeleton vertex is strictly positive and is the same for all $x\in {\mathbb Z}$.
\item\vskip2pt
The sequence of events 
$\{x\in {\mathcal S}\}$, $x=\ldots,-2,-1,0,1,2,\ldots\strut$ is stationary ergodic and, therefore,
with probability one there are infinitely many skeleton points $\{t_i\}$.
\item\vskip2pt
The sequence
$\{t_i\}$ (where $\ldots t_{-2}<t_{-1}<0\le t_0< t_1<\ldots\strut$) forms a stationary renewal process (in discrete time) and, in particular, the lengths of the intervals $\{t_i-t_{i-1},\:
i\in\mathbb{Z},\: i\ne 0\}$ are independent and identically distributed random variables that do not depend on the pair of random variables $(t_{-1},t_0)$, and the latter random variables depend on each other and $t_0$ has the same distribution as $|t_{-1}|-1$. Further, 
$\P(t_{-1}=-i)=\P(t_1-t_0\ge i)/\E(t_1-t_0)$, $i=1,2,\ldots\strut$.
\item\vskip2pt
For some 
$C>0$,
\begin{equation}\label{09.01.1}
\E e^{Ct_0}<\infty\quad \text{and hence}\quad \E e^{C(t_1-t_0)}<\infty.
\end{equation}
\item
For any $j<k$ denote by 
$$
L_{j,k}=\max_{\pi\in\Pi_{j,k}} L(\pi)\quad (\text{where the maximum over an empty set
is}\ -\infty)
$$
the maximal path length among the paths from 
$\Pi_{j,k}$. Then, for each $n>0$, if $L_{0,n}>0$, then any path of length  $L_{0,n}$ from $0$ to $n$ has to include all intermediate skeleton points (if there are any). Namely, assume that $0\le t_0<t_1< \ldots<t_m\le
n<t_{m+1}$, for some $m\ge 0$. Then any path of maximal length from $0$ to $n$, that belongs to the set $\Pi_{0,n}$,  has to include every vertex $t_0,\ldots,t_m$. Further, with necessity, all values 
$L_{t_0,t_1},\ldots,L_{t_{m-1},t_m}$ are strictly positive and 
$$
L_{0,n}=L_{0,t_0}+L_{t_0,t_1}+\ldots+L_{t_{m-1},t_m}+L_{t_m,n}.
$$
\end{enumerate}

Along with the set of paths $\Pi_{j,k}$ (introduced in Section~1) that include only edges of finite weights, we introduce also the set of paths $\Pi^+_{j,k}$ from $j$ to $k$ that include only edges with positive weights (i.e. $v_e>0$ for every $e\in\pi$, given
$\pi\in\Pi^+_{j,k}$) and let 
$$
w_{j,k}^+:=\max\limits_{\pi\in\Pi_{j,k}^+} w(\pi).
$$

\begin{definition}
Vertex $x\in\mathbb{Z}$ is a \emph{skeleton-plus} vertex if it is connected to any other vertex by a path with edges of positive weights only, i.e., for and $j<x$ and $k>x$, there is a path $\pi$ from $j$ to $x$ and a path $\widetilde{\pi}$ from  $x$ to $k$ such that
$v_e>0$ for all $e\in\pi$ and all $e\in\widetilde{\pi}$ and, in particular, inequalities $w_{j,x}^+> 0$ and $w_{x,k}^+> 0$ hold. Denote by
$\mathcal{S}^+=\{t_i^+\}$ the set of skeleton-plus vertices.
\end{definition}

Notice that since we assume that $p^+>0$, the results from the paper \cite{Foss2} are also applicable to the set $\mathcal{S}^+$ and, in particular, exponential moments in \eqref{09.01.1} stay finite if one replaces $\{t_i\}$ by $\{t_i^+\}$.

Now we introduce the sets 
$\mathcal{R}$ of \emph{renewal} and $\mathcal{R}^+$ of \emph{renewal-plus} vertices. Let $c_1\ge c_2>0$ be two fixed numbers. For 
$x\in\mathbb{Z}$, introduce the following events:
\begin{align*}
A_x^r(c_1)&=\bigcap\limits_{i=1}^\infty\{w_{x,x+i}\ge c_1i\},\\
A_x^{0}(c_2)&=\bigcap\limits_{j,i=1}^\infty\{v_{x-j,x+i}< c_2(j+i)\},\\
A_x^l(c_1)&=\bigcap\limits_{j=1}^\infty\{w_{x-j,x}\ge c_1j\}.
\end{align*}

Here event $A_x^r(c_1)$ means that vertex $x$ is connected to any vertex on its right by a path of a finite weight and, moreover, all corresponding maximal path weights are strictly positive and increase at least linearly ( with speed not smaller than $c_1$) with growth of distance from $x$.
Similarly, event $A_x^l(c_1)$ means that vertex $x$ is connected to any vertex on its left by a path of a finite weight and, moreover, all corresponding maximal path weights are strictly positive and grow at least linearly with growth of distance from $x$. In particular, if both events $A_x^r(c_1)$ and $A_x^l(c_1)$ occur, then vertex 
$x$ is necessarily a skeleton vertex and, in addition, $w_{x-j,x+i}\ge c_1(j+i)$, for all $j,i\ge 0$. If, in addition, event $A_x^0(c_2)$ takes place, then (by $c_2\le c_1$), for any $j,i\ge 0$,
any path of maximal weight from vertex 
$x-j$ to vertex $x+i$ must necessarily include vertex~$x$.

\begin{definition}
Vertex $x$  is a \emph{renewal} vertex if all three events $A_x^l(c_2)$, $A_x^0(c_1)$ and  $A_x^r(c_2)$ occur. Let  $\mathcal{R}=\{\Gamma_i\}$ be the random set of renewal vertices.
\end{definition}

Since the events $\{x\in\mathcal{R}\}$ form a stationary ergodic sequence, the ``0-1'' law holds: with probability one, the set~$\mathcal{R}$ is either infinite or empty. Assume it is infinite, then its elements may be ordered as
$$
\ldots<\Gamma_{-1} <0\le\Gamma_0<\Gamma_1<\ldots,
$$
and representation \eqref{mainpres} holds. It is this representation that helps us to study the asymptotics of sequence $w_{0,n}$, as $n$ grows to infinity.

In analogy to the notation above, we introduce events
\begin{align*}
A_x^{r+}(c_1)&=\bigcap\limits_{i=1}^\infty\{w_{x,x+i}^+\ge c_1i\},\\
A_x^{0+}(c_2)\equiv A_x^{0}(c_2)&=\bigcap\limits_{j,i=1}^\infty\{v_{x-j,x+i}<
c_2(j+i)\},\\ A_x^{l+}(c_1)&=\bigcap\limits_{j=1}^\infty\{w_{x-j,x}^+\ge c_1j\}.
\end{align*}
Here the event $A_x^{r+}(c_1)$ means that vertex $x$ is connected to any vertex on its right by a path with edges of positive weights only and, moreover, all related maximal path weights grow at least linearly with the distance from vertex~$x$. Similarly, the event $A_x^{l+}(c_1)$ means that vertex $x$
is connected to any vertex on its left by  a path with edges of positive weights only and, moreover, all related maximal path weights grow at least linearly with the distance from vertex~$x$.  Necessarily, $A_x^{l+}(c_1) \subseteq
A_x^l(c_1)$ and $A_x^{r+}(c_1)\subseteq A_x^r(c_1)$, for any $c_1>0$.

\begin{definition} Call $x$ a \emph{renewal-plus} vertex if
$A_x^{l+}(c_1)\cap A_x^{0+}(c_2)\cap A_x^{r+}(c_1)$ occurs, and denote by $\mathcal{R}^+$ 
the set of all renewal-plus vertices. Note that this set is also either infinite a.s. or empty a.s.
\end{definition}

The following relations hold:
\begin{equation}\label{relations}
\mathcal{R}^+ \subseteq \mathcal{R} \subseteq \mathcal{S}\qquad \text{and}\qquad
\mathcal{R}^+ \subseteq \mathcal{S}^+ \subseteq \mathcal{S}.
\end{equation}
In addition, the sets $\mathcal{R}$ and  $\mathcal{R}^+$ increase as $c_2$ increases and  $c_1$ decreases.

\subsection{Regenerative structure and existence of exponential moments}

For formulation of the next statement, we need the distribution of random variable $v^+$ that was defined in \eqref{vplus}. Let $\essinf v^+=\inf \{t>0 :\: \P(v^+<t)>0\}$ and 
$V=\E\min\limits_{t_0^+\le i<j\le t_1^+} v_{i,j}^+$. Clearly, $\essinf v^+<V$ if the distribution of random variable $v^+$ does not degenerate. Let 
$\gamma^+=\nfrac{1}{\E(t_1^+-t_0^+)}$. From Lemmas 5--7 in \cite{Foss2}, we get the following result.

\begin{lemma}\label{l1}
Assume that condition \eqref{maincond} holds. If 
\begin{equation}\label{defc}
\gamma^+ \essinf v^+<c_2\le c_1<\gamma^+ V,
\end{equation}
then, for any
$x\in\mathbb{Z}$,  $\P(A_x^{l+}(c_1)\cap A_x^{0}(c_2)\cap A_x^{r+}(c_1))>0$  and  the set 
$\mathcal{R}^+$ is infinite with probability $1$. Therefore, the set $\mathcal{R}$ is infinite with probability $1$, too.
\end{lemma}

\begin{remark}
In paper \cite{Foss2}, renewal points have been introduced in the case $c_1=c_2$ and corresponding results have been proved in this case only. However, one can easily check that the proofs of these results remain unchanged (apart from minor changes in notation) under the more general conditions \eqref{defc}.
\end{remark}

Introduce the following cycles
$$
\mathcal{C}^+_k:=\Bigl(\Gamma^+_k-\Gamma^+_{k-1};\:
\Bigl\{v^+_{\Gamma^+_{{k-1}+j},\Gamma^+_{{k-1}+i}},\: 0\le
j<i\le\Gamma^+_k-\Gamma^+_{k-1}\Bigr\}\Bigr),\quad k\in\mathbb{Z},
$$
and
$$
\mathcal{C}_k:=\bigl(\Gamma_k-\Gamma_{k-1};\:
\bigl\{v_{\Gamma_{{k-1}+j},\Gamma_{{k-1}+i}},\: 0\le
j<i\le\Gamma_k-\Gamma_{k-1}\bigr\}\bigr),\quad k\in\mathbb{Z}.
$$

The following results take place.

\begin{lemma}\label{l2} Under conditions 
\eqref{maincond}  and \eqref{defc}, the following two statements hold:
\begin{enumerate}[\rm(I)]
\item\vskip-5pt
Random elements 
 $\{\mathcal{C}^+_k,\: k\in\mathbb{Z}\}$ are mutually independent and random elements $\{\mathcal{C}^+_k,\:
k\in\mathbb{Z}\setminus\{0\}\}$ are identically distributed. The process
$\Bigl(\Gamma^+_k,w^+_{\Gamma^+_{k-1},\Gamma^+_k}\Bigr)$
$k\in\mathbb{Z}$ is a stationary marked point process in discrete time and it generates a stationary CRP, i.e. its first coordinates $\Gamma^+_k$ form a stationary point process with corresponding marks
$w^+_{\Gamma^+_{k-1},\Gamma^+_k}$.
\item
The previous statement remains valid for cycles
 $\mathcal{C}_k$, with natural replacements of $\Gamma^+_k$ by $\Gamma_k$ and 
$w^+_{\Gamma^+_{k-1},\Gamma^+_k}$ by $w_{\Gamma_{k-1},\Gamma_k}$.
\end{enumerate}
\end{lemma}

Note that the first statement of Lemma \ref{l2} is a direct consequence of Lemma
3.8 from \cite{Foss1} and that the schemes of the proofs of both statements of Lemma \ref{l2}
are identical.

We will formulate and prove now the first of the main statements of this Section.

\begin{lemma}\label{ll3}
Assume that conditions \eqref{maincond} and \eqref{defc} hold. Then there exists constant $C>0$ 
such that
\begin{equation}\label{gamma1+}
\E e^{C\Gamma^+_0}<\infty\quad \text{and hence}\quad \E
e^{C(\Gamma^+_1-\Gamma^+_0)}<\infty.
\end{equation}
Since the sequence 
$\{\Gamma_n^+\}$ is a subsequence of $\{\Gamma_n\}$, the statement  \eqref{gamma1+} remains valid after the replacement of $\Gamma_0^+$ and
$\Gamma_1^+$ by, correspondingly, $\Gamma_0$ and $\Gamma_1$, i.e. both first expectations in 
\eqref{tauzeta2} are finite, for some $C>0$.
\end{lemma}

\begin{proof} We borrow from 
\cite{Foss1} a number of auxiliary constructions. Introduce a set 
$\mathcal{U}=\{x\in\mathbb{Z}:\: \mathbf{I}(A_x^{l+}(c_1))=1\}$. It is not difficult to see that 
$\mathcal{R}^+\subseteq\mathcal{U}$. Enumerate the elements of the set $\mathcal{U}$ in the increasing order $\ldots,\rho_{-1},\rho_0,\rho_1,\ldots\strut$ 
%$\mathcal{U}$, 
where $\rho_0$ is its smallest non-negative element. Define a new sequence of cycles: for $k\in\mathbb{Z}$,
$$
\mathcal{D}_k=\bigl(\rho_k-\rho_{k-1};\: \bigl\{v^+_{\rho_{k-1}+j,\rho_{k-1}+i},\:
0\le j<i\le\rho_k-\rho_{k-1}\bigr\}\bigr).
$$

The next result folllows from Lemma 3.10 of \cite{Foss1}.

\begin{lemma}\label{l3}
Assume that conditions \eqref{maincond} and \eqref{defc} hold. Then cycles 
$(\mathcal{D}_k,\: k\in\mathbb{Z})$ are mutually independent random elements where   
$(\mathcal{D}_k,\: k\in\mathbb{Z}\setminus\{0\})$ are identically distributed, and  sequence $\rho_n$,
$n\in\mathbb{Z}$ forms a stationary point process that generated the corresponding CRP.
\end{lemma}
The statement of Lemma 
\ref{l3} may be easily explained in simple terms. For that, introduce for $d>0$ events
$$
A_{x,d}^{l+}(c_1):=\bigcap\limits_{j=1}^d\bigl\{w^+_{x-j,x}\ge c_1j\bigr\}.
$$
Note that, for any integers $k\ge 0$ and $0\le s_0< s_1<\ldots<s_k$, the event $\{\rho_0=s_0,\ldots,\rho_k=s_k\}$ is uniquely determined by the collection of random variables
$\mathcal{B}_{s_k}:=\{v_{i,j},\: i<j\le s_k\}$ and, on this event,
the equality $\rho_{k+1}-\rho_k=m$ holds if and only if $m=\min
\bigl\{d>0:\: \mathbf{I}\bigl(A^{l+}_{s_k,s_k+d}(c_1)\bigr)=1\bigr\}$, and the latter is determined by the random variables 
$\mathcal{B}_{s_k,s_k+m}:=\{v^+_{i,j},\: s_k\le i<j\le m\}$. It is not difficult to see that the families of random variables $\mathcal{B}_{s_k}$ and 
$\mathcal{B}_{s_k,s_k+m}$ are mutually independent and that the distribution of random variables from $\mathcal{B}_{s_k,s_k+m}$ does not depend on $k$ and $s_k$. Essentially these facts imply the statement of Lemma \ref{l3}.

For $d>0$, introduce auxiliary events
\begin{align*}
A_{x,d}^{r+}(c_1) &:=\bigcap\limits_{i=1}^d\{w^+_{x,x+i}\ge c_1i\},\\
A_{x,d}^{0+}(c_2) &:=\bigcap\limits_{1\le i\le d,\, j\ge 1}^\infty\{v_{x-j,x+i}<
c_2(j+i)\}.
\end{align*}

Let $\mu =\inf\bigl\{d>0:\:\mathbf{I}\bigl(A_{0,d}^{0+}(c_2)\cap
A_{0,d}^{r+}(c_1)\bigr)=0\bigr\}$. Note that $\P(\mu=\infty)= \P(A_0^{0+}(c_2)\cap
A_0^{r+}(c_1))>0$. Define recursively random variables $\sigma_0,\mu_0, \ldots, \sigma_K,\mu_K,$ where $K=\min \{k\ge 0:\:
\mu_k=\infty\}$. Let 
$$
\sigma_0=\rho_0,\qquad
\mu_0=\inf\Bigl\{d>0:\:\mathbf{I}\bigl(A_{\sigma_0,\sigma_0+d}^{0+}(c_2)\cap
A_{\sigma_0,\sigma_0+d}^{r+}(c_1)\bigr)=0\Bigr\},
$$
and, for
$k=0,1,\ldots\strut$, if $\mu_k<\infty$, then 
\begin{align*}
\sigma_{k+1}&=\inf\{x\in\mathcal{U}:\: x\ge\sigma_k+\mu_k\},\\*
\mu_{k+1}&=\inf\Bigl\{d>0:\: \mathbf{I}\bigl(A_{\sigma_k,\sigma_k+d}^{0+}(c_2)\cap
A_{\sigma_k,\sigma_k+d}^{r+}(c_1)\bigr)=0\Bigr\}.
\end{align*}

The process of construction of this sequence is presented more transparently in Figure~\ref{fig1}.

\begin{figure}[tp]
\centering\vskip3pt
\setlength{\unitlength}{2092sp}%
\begin{picture}(0,0)(1489,-4574)
\put(2296,-4486){\makebox(0,0)[lb]{$\rho_0$}}
\put(1741,-4516){\makebox(0,0)[lb]{$0$}} \put(2866,-3556){\makebox(0,0)[lb]{$\mu_0$}}
\put(5581,-4541){\makebox(0,0)[l]{$\underset{\textstyle\sigma_1}{\rho_1}$}}
\put(6196,-3421){\makebox(0,0)[lb]{$\,\,\,\mu_1$}}
\put(6481,-4501){\makebox(0,0)[lb]{$\rho_2$}}
\put(7966,-4541){\makebox(0,0)[l]{$\underset{\textstyle\sigma_2}{\rho_3}$}}
\put(8371,-3721){\makebox(0,0)[lb]{$\,\mu_2$}}
\put(9166,-4541){\makebox(0,0)[l]{$\underset{\textstyle\sigma_3}{\rho_4}$}}
\put(10096,-4501){\makebox(0,0)[lb]{$\rho_5$}}
\put(11101,-2911){\makebox(0,0)[lb]{$\mu_3=\infty$}}
\end{picture}%
\includegraphics[width=.85\textwidth]{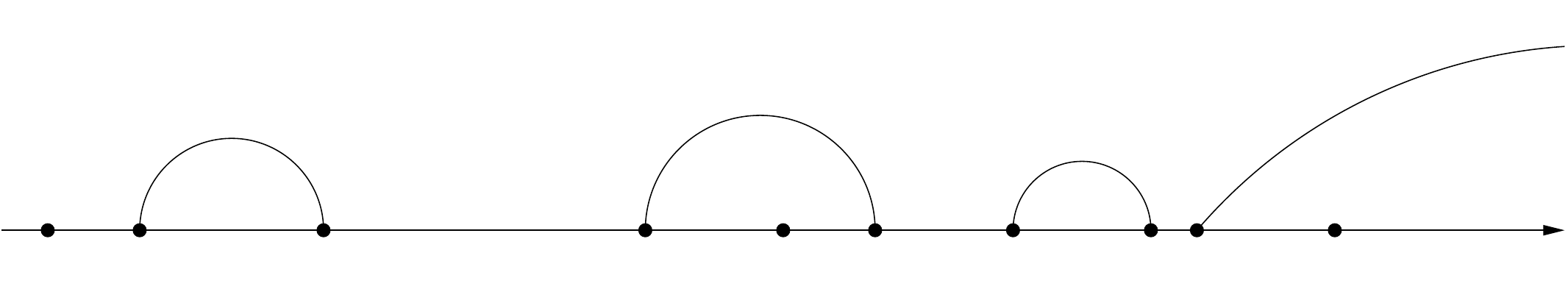}\vskip3pt
\caption{The process of construction of sequence $\sigma_k$. Here $K=3$ and $\sigma_3=\rho_4$ is the moment of regeneration that goves an upper bound for $\Gamma_0$}.\label{fig1}
\end{figure}
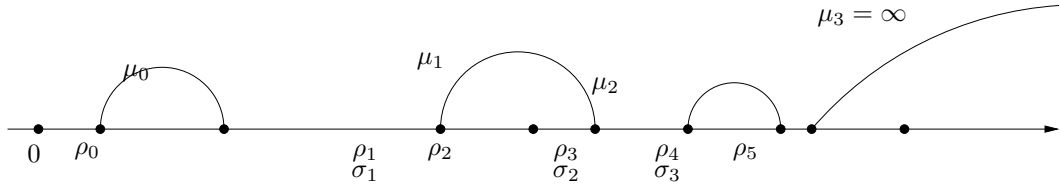

As it follows from the construction, $\sigma_k\in\mathcal{U}$ for $k\le K$ and $\sigma_K\in\mathcal{R}^+$ and, therefore, 
$\sigma_K\ge\Gamma_0^+$ a.s. Moreover, the random variables $\{\mu_k\}$ form a sequence of independent and identically distributed random variables having a common distribution with $\mu$. 
Therefore, the random variable $K$ is geometrically distributed with parameter $q$ and, in particular,
its exponential moments $\E e^{CK}$ are finite for $C<-1/\ln
q$. In addition, for any $k\ge 1$, given $\{K=k\}$, the random variables $\mu_0, \ldots, \mu_{k-1}$ are conditionally independent and distributed as $\P(\mu_0\in\cdot\mmid
\mu_0<\infty)$.

Since $\rho_k-\rho_{k-1}\ge 1$,
\begin{equation}\label{11.01.2}
\sigma_{K}\le\rho_M=\rho_0+\sum\limits_{k=1}^M(\rho_k-\rho_{k-1}),
\end{equation}
where $M:=\sum\limits_{j=0}^{K-1}\mu_j$.

Therefore, if we show that the term on the right-hand side of \eqref{11.01.2} has a finite exponential moment, then $\sigma_K$ also has a finite exponential moment, and this in turn implies the statement of Lemma \ref{ll3}. Given Lemma \ref{l4} from Section~4 and the elementary inequality
$e^{x+y}< e^{2x}+e^{2y}$
for $x,y\ge 0$, it is sufficient to show that, firstly, random variable $\rho_1-\rho_0$ has also a finite exponential moment (then the same holds for $\rho_0$) and, secondly, the probabilities $\P(\mu_0=m\mmid \mu_0<\infty)$ decay exponentially fast in $m$.

Due to the independence of $\rho_0$ and $\rho_1-\rho_0$ and since, for any natural $n>0$, given occurrence of event $A_0^{l+}(c_1)$, events $A_n^{l+}(c_1)$ and $A_{n,n}^{l+}(c_1)$ either occur or not simultaneously, we get that the distributions of random variables 
$\rho_1-\rho_0$ and $\nu := \min\{n:\: \mathbf{I} (A_{n,n}^{l+}(c_1))=1\}$ coincide:
$$
\P(\rho_1-\rho_0=m)=\P(\rho_1-\rho_0=m\mmid \rho_0=0)=\P(\nu=m\mmid \rho_0=0)=\P(\nu
=m),\quad m=1,2,\ldots,
$$
and existence of an exponential moment of random variable 
$\nu$ follows from Proposition 3.12 in \cite{Foss1}.

Further,
\begin{align*}
\P(\mu=d)&=\P\Bigl(\bigl(A_{0,d}^{0+}(c_2)\cap
A_{0,d}^{r+}(c_1)\bigr)^c\cap\bigl(A_{0,d-1}^{0+}(c_2)\cap
A_{0,d-1}^{r+}(c_1)\bigr)\Bigr)\\ &\le\P\Bigl(\bigl(A_{0,d}^{0+}(c_2)\bigr)^c\cap
A_{0,d-1}^{0+}(c_2)\Bigr)+ \P\Bigl(\bigl(A_{0,d}^{r+}(c_1)\bigr)^c\cap
A_{0,d-1}^{r+}(c_1)\Bigr)\\ &\le\P\Bigl(\,\sup\limits_{j\ge
1}\bigl(v^+_{-j,d}-c_2j\bigr)>c_2d\Bigr)+ \P(w_{0,d}^+<c_1d)\\[-2pt] &\le
\sum\limits_{j=1}^\infty\P\bigl(v>c_2(d+j)\bigr)+ \P\bigl(w^+_{0,d}<c_1d\bigr).
\end{align*}
In the last line, the sum of probabilities 
$\sum\limits_{j=1}^\infty\P(v>c_2(d+j))$ decays in $d$ exponentially fast, thanks to  \eqref{maincond}. In order to show that the very last probability decreases exponentially fast in $d$ too, 
we choose $\varepsilon >0$ such that 
$\widetilde{c}:= c_1(1+3\varepsilon)$ also satisfies \eqref{defc}. Let $\eta (d)=\max \{k:\: t_k^+\le d
\}$ (where, by convention, the maximum over empty set is $-\infty$). Then, for $r= \gamma^+ (1+\varepsilon)^{-1}$,
\begin{align*}
\P\bigl(w^+_{0,d}<c_1d\bigr)&\le\P (t_0^+>d)+\P\Bigl(w^+_{t_0^+,
t^+_{\eta(d)}\!}<c_1d,\: t_0^+\le d\Bigr)\\ &\le\P (t_0^+>d) +
\P\bigl(\eta(d)<[rd]\bigr) + \P\Biggl(\,\sum_{i=1}^{[rd]}
w^+_{t_{i-1},t_i}<c_1d\Biggr)\\ &\le\P (t_0^+>d)+\P\Biggl(\,\sum_1^{[rd]}
(t_i^+-t_{i-1}^+)> d\Biggr)+\P\Biggl(\,\sum_{i=1}^{[rd]}
w^+_{t_{i-1},t_i}<c_1d\Biggr),
\end{align*}
where $[rd]$ is the integer part of $rd$. In the last line of these inequalities all three summands decrease exponentially fast in $d$: the first summand since $t_0^+$ has a finite exponential moment, 
the second summand since the increments $t_i^+-t_{i-1}^+$ have a finite exponential moment, $\E(t_1^+-t_0^+) r<1$, and, by the exponential Chebyshev inequality with $h>0$, 
$$
\P\Biggl(\,\sum_{i=1}^{[rd]} (t_i^+-t_{i-1}^+)> d\Biggr)\le\Bigl(\bigl(\E\exp
\bigl(h(t_1^+-t_0^+)\bigr)\bigr)^re^{-h}\Bigr)^d,
$$
where the right-hand side decays exponentially fast in $d$, if one takes~$h>0$ sufficiently small.
Finally, the third summand decays exponentially fast because of a well-known fact:
for any sequence $X,X_1,X_2,\ldots\strut$ of independent and identically distributed positive random variables with finite mean $\E X$ and for any $\delta\in (0,1)$ the probabilities 
$\P\Bigl(\,\sum\limits_{i=1}^n X_i<(1-\delta)n\E X\Bigr)$  decrease exponentially fast as $n$ grows. In our case, $n=[rd]\ge rd-1$,
$X_i=w_{t_{i-1}^0,t_i^0}$, $\E w_{t_0^+,t_1^+}\ge V$ and $c_1d\le
c_1(1+n)(1+\varepsilon)/\gamma^+\le c_1(1+2\varepsilon)n V/\gamma^+<(1-\delta)nV$ for all sufficiently large
$n$, where $\delta=\varepsilon/(1+3\varepsilon)$.

Thus, the probabilities $\P(\mu =d)$ and, therefore, the probabilities $\P(\mu =d\mmid \mu<\infty)$ decrease exponentially fast as $d$ grows. This completes the proof of Lemma \ref{ll3}.\qed
\end{proof}

We will proceed now with the proof of finiteness of the two last mathematical expectations 
in \eqref{tauzeta2}.

\begin{lemma}\label{wfinite}
Assume that conditions \eqref{maincond} and \eqref{defc} hold. Then $\E\exp(Cw_{0,\Gamma_0})<\infty$ and, therefore, 
$\E\exp (Cw_{\Gamma_0,\Gamma_1})<\infty$, for some $C>0$.
\end{lemma}

\begin{proof} 
Choose any path 
$\pi$ from vertex $0$ to vertex $\Gamma_0$ and assume that it includes $d+1$ vertices, $0=x_0<x_1<\ldots<x_d=\Gamma_0$.
Since $\sum\limits_{k=1}^d(x_k-x_{k-1})=\Gamma_0$, we get 
\begin{equation}\label{11.01.7}
\begin{aligned}[b]
w_{0,\Gamma_0}=\sum\limits_{k=1}^dv_{x_{k-1},x_k}
&\le\Gamma_0+\sum\limits_{k=1}^d(v_{x_{k-1},x_k}-(x_k-x_{k-1}))^+\\
&\le\Gamma_0+\sum\limits_{0\le x <y\le\Gamma_0}(v_{x,y}-(y-x))^+\\
&\le\Gamma_0+\sum\limits_{x=0}^{\Gamma_0-1}Z_x,
\end{aligned}
\end{equation}
where $\bigl\{Z_x:=\max\limits_{y>x}(v_{x,y}-(y-x))^+\bigr\}_{x\in\mathbb{Z}}$ is a sequence
of independent and identically distributed random variables. By condition 
\eqref{maincond}, the tail distribution 
$$
\P(Z_0>m)\le\sum_{k=1}^{\infty} \P(v>m+k)
$$
decays exponentially fast in 
$m$. To complete the proof, it is enough to use the inequality $e^{x+y}<e^{2x}+e^{2y}$  and Lemma \ref{l4}.\qed
\end{proof}

We complete Section~2 with a short proof of a simple fact. 

\begin{lemma}\label{l9}
Let $p\in(0,1]$ and assume that conditions \eqref{maincond} and \eqref{defc} hold. Then 
$\P(\Gamma_1-\Gamma_0=1,\: w_{\Gamma_0,\Gamma_1}\ge y)>0$ for any $y\in(c_2,\esssup v^+)$.
\end{lemma}

\begin{proof}
The following two events coincide:
\begin{equation}\label{5events}
\{\Gamma_0=0,\: \Gamma_1-\Gamma_0=1,\: w_{\Gamma_0,\Gamma_1}\ge y\} =A^l_0(c_1) \cap
A^{0+}_{0,1}(c_2) \cap \{v_{0,1}\ge y\} \cap B_{1,1}(c_2) \cap A^r_1(c_1),
\end{equation}
where 
$$
B_{1,1}(c_2)=\bigcap_{j=1}^{\infty} \{v_{0,1+j}<c_2(1+j)\}.
$$
Since all five events in the right-hald side of 
\eqref{5events} are mutually independent and each of them has positive probability, the result follows.\qed
\end{proof}

\section{Analysis of the compound renewal process and proof of the main theorem}

We analyse here the CRP determined by the stationary marked process
$(\Gamma_k,w_{\Gamma_{k-1},\Gamma_k})$, $k\in\mathbb{Z}$. Using results from the previous Section and from \cite{Mog1}, and also the classical Stone's theorem \cite{St}, we will show that, for any admissible constants $c_1,c_2$, the corresponding  CRP has the same exact asymptotics with sequence  $w_{0,n}$ in the domain of the normal and moderately large deviations (this is the result of Theorem  1). Next, we apply a corresponding change of measure in order to remove the ``defect'' of the CRP. Then we conclude with the statement that parameters $\alpha$, $\sigma^2$ and $D(\alpha)$ that appear in Theorem 1, in fact, do not depend on constants $c_1$ and $c_2$.

In what follows, it will be convenient to us to introduce some notation that in corrrespondence with notation from \cite{Mog1}:
\begin{equation}\label{1.1}
(\tau_k, \bfm{u}_k):=(\tau_k,(u_{k,1},\ldots,u_{k,\tau_k})),\quad k=1,2,\ldots,
\end{equation}
where
\begin{align*}
(\tau_1,(u_{1,1},\ldots,u_{1,\tau_1}))&:=
(\Gamma_0,(w_{0,1},\ldots,w_{0,\Gamma_0})),\\
(\tau_k,(u_{k,1},\ldots,u_{k,\tau_k}))&:=
(\Gamma_{k-1}-\Gamma_{k-2},(w_{\Gamma_{k-2},\Gamma_{k-2}+1},
\ldots,w_{\Gamma_{k-2},\Gamma_{k-1}})),\quad \text{for}\ k\ge 2.
\end{align*}
It was shown in Section~2 that vectors $(\tau_k, \bfm{u}_k)$, $k\ge 2$, are indentically distributed, and we will use notation 
\begin{equation}\label{03.02.2}
(\tau, \bfm{u}):=(\tau,~(u_1,\ldots,u_{\tau}))
\end{equation}
for any vector having this distribution. 
We let also $\zeta=u_{\tau}$ and $\zeta_k=u_{\tau_k}$, for $k\ge 1$. Then, in particular, $\{(\tau_k,\zeta_k)\}$ is a sequence of independent random vectors that have the same distribution for  $k\ge
2$ with vector $(\tau,\zeta)$.

We will list now statements from 
Section~2 (based on Lemmas \ref{l2}, \ref{ll3}, \ref{wfinite} and \ref{l9}) that we need here.
Let $p\in (0,1]$ and let $c_1,c_2$ satisfy condition 
~\eqref{defc}.
\begin{enumerate}[$S_{\mathrm{III}}$.]
\item[$S_{\mathrm{I}}$.]\vskip-5pt
The sequence 
\eqref{1.1} consists of independent random vectors, and,for $k\ge 2$, the random vectors $(\tau_k,\bfm{u}_k)$ have the same distribution.
\item[$S_{\mathrm{II}}$.]
The random variables
$u_{1,\tau_1}$ and $u_1,\ldots,u_{\tau}$ are positive, and
$$
\max\{u_1,\ldots,u_{\tau-1}\}\le u_{\tau}.
$$
\item[$S_{\mathrm{III}}$.]
For some 
$C>0$,
$$
\E e^{C\tau_1}<\infty,\quad \E e^{C\tau}<\infty,\quad \E
e^{Cu_{1,\tau_1}}<\infty,\quad \E e^{Cu_{\tau}}<\infty.
$$
\item[$S_{\mathrm{IV}}$.]
The probability 
$\P(\tau=1,u_{\tau}\ge y)$ is strictly positive for any $y\in(c_2,\esssup v^+)$ .
\end{enumerate}\vskip-5pt
Denote 
\begin{equation}\label{1.4}
(\tau,\zeta):=(\tau,u_{\tau}),\qquad (\tau_k,\zeta_k):=(\tau_k,u_{k,\tau_k}),\quad
k=1,2,\ldots%.
\end{equation}
Then the sequence 
$\{(\tau_k,\zeta_k)\}$ consists of independent random vectors that have, for $k\ge 2$,
a common distribution with vector $(\tau,\zeta)$.
The results above imply  
\begin{corollary}\label{r1}
Let conditions 
\eqref{maincond} and \eqref{defc} hold. Then the following statements take place. 
\begin{enumerate}[\rm(I)]
\item\vskip-5pt
If the random variable $v$ satisfies condition $[\mathbf{Z}]$, then 
\begin{itemize}
\item[$\lbrack{\bf ZZ}\rbrack$]
The distribution of the random vector $(\tau,\zeta)$ is arithmetic and is concentrated on the lattice of span $1$ for each of its coordinates.
\end{itemize}
\item
If the random variable $v$ satisfies condition $[{\bf R}]$, then $(\tau,\zeta)$
\begin{itemize}
\item[$\lbrack{\bf ZR}\rbrack$]
The marginal distribution of the first coordinate of vector $(\tau,\zeta)$ is arithmetic with span 1
and the marginal distribution on the second coordinate is non-lattice.\footnote{In terms of the characteristic function $f(z,l):=\E e^{iz\tau+il\zeta}$ of random vector  $(\tau,\zeta)$, these two conditions may be presented as:\\
$[{\bf ZZ}\,]$ $f(2\pi z,2\pi t)=1$ for any $(z,t)\in\mathbb{Z}^2$ and $|f(z,t)|<1$ for any 
$(z,t)\notin\mathbb{Z}^2$.\\
$[{\bf ZR}]$ $f(2\pi z,0)=1$ for any $z\in\mathbb{Z}$, $|f(2\pi
z,0)|<1$ for any $z\notin\mathbb{Z}$ and $|f(0,t)|<1$ for any $t\ne 0$.}.
\end{itemize}
\end{enumerate}
\end{corollary}

We turn now to the proof of the main result.

\proofof{Theorem~\ref{th2}}
Consider sequence 
$\left\{(\tau_k,\zeta_k)\right\}_{k=1}^\infty$ of independent random vectors, having for $k\ge 2$ 
the same distribution as $(\tau,\zeta)$ (see formula \eqref{03.02.2} and notation after it).
Introduce sequences of partial sums 
$$
T_n:=\sum_{k=0}^n\tau_k,\qquad Z_n:=\sum_{k=0}^n\zeta_k,\quad n\ge 0,
$$
where $(\tau_0,\zeta_0):=(0,0)$. Let 
\begin{align*}
\eta_+(n)&:=\min\{k\ge 1:\: T_k>n\},\\ \nu_+(n)&:=\max\{k\ge 0:\: T_k\le
n\}=\eta_+(n)-1,\\ \gamma_+(n)&:=n-\nu_+(n).
\end{align*}
Introduce our CRP (we will call it the ``first CRP'') by 
$$
Z_+(n):=\sum_{k=0}^{\nu_+(n)}\zeta_k.
$$
Let us clarify that we consider here notation $Z_+(n)$ with low case $+$ just to reconcile our notation with that from the paper~\cite{Mog1}, where, in addition to notation $Z_+(n)$, $\nu_+(n)$, $\gamma_+(n)$, further notation $Z(n)$, $\nu(n)$,
$\gamma(n)$ have been used. We like to point out that thus ``plus'' does not relate to that in notation  $w_{j,m}^+$.

Using the introduced notation, we obtain the representation 
$$
w_{0,n}=Z_+(n)+w_{n-\gamma_+(n),n}=Z_{\nu_+(n)}+w_{\nu_+(n),n}.
$$
Consider a random vector (of random length), given by formula \eqref{03.02.2}:
$$
(\tau,(u_1,u_2,\cdots,u_\tau)),
$$
whose coordinates clearly satisfy the following constraints:
$$
\tau\ge 1,\qquad \min_{1\le i\le\tau}u_i\ge c_1>0,\qquad \max_{1\le i\le\tau}u_i\le
u_\tau.
$$
Next, introduce a random vector 
$(\tau^*,\zeta^*)$, taking values  $(i,y)\in\mathbb{Z}\times \mathbb{R}$, $i\ge 0$, $y\ge 0$, by 
$$
\P(\tau^*=i,\: \zeta^*\in dy):=\frac{1}{Q} \P(\tau\ge i+1,\: u_i\in dy),
$$
(we assume 
$u_0=0$ a.s.), where
$$
Q:=\sum_{i=0}^\infty\,\int_0^\infty\P(\tau\ge i+1,\: u_i\in dy)=
\sum_{i=0}^\infty\P(\tau\ge i+1)=\mathbf{E}\tau.
$$
It follows from condition
$S_{\mathrm{III}}$ that there exists a constant $C>0$ such that 
$$
\E e^{C\tau^*}<\infty,\qquad \E e^{C\zeta^*}<\infty.
$$
Along with the sequence
$\{(\tau_k,\zeta_k)\}$ that determines CRP $Z_+(n)$ and the functionals $\nu_+(n)$, $\gamma_+(n)$,
we introduce another sequence $\{(\tau^*_k,\zeta^*_k)\}$ by 
$$
(\tau^*_1,\zeta^*_1):=(\tau_1,\zeta_1)+(\tau^*,\zeta^*),\qquad
(\tau^*_k,\zeta^*_k):=(\tau_k,\zeta_k)\quad \text{for}\ k\ge 2,
$$
where the random vector 
$(\tau^*,\zeta^*)$ does not depend on the sequence $\{(\tau_k,\zeta_k)\}$. The new sequence $\{(\tau^*_k,\zeta^*_k)\}$ determines a new CRP $Z_+^*(n)$ and new functionals 
$\nu_+^*(n)$ and~$\gamma_+^*(n)$.

\begin{lemma}\label{l10} 
Assume that conditions $S_{\mathrm{I}}$--$S_{\mathrm{IV}}$ hold. Then, for any $n\ge 2$ and any real $x\ge c_1$ and $\Delta>0$, the following equality takes place:
\begin{equation}\label{1.6}
\P\bigl(Z_+(n)+w_{n-\gamma_+(n),n}\in [x,x+\Delta),\: \tau_1\le n\bigr)=
Q\P\bigl(Z_+^*(n)\in [x,x+\Delta),\: \gamma_+^*(n)=0\bigr).
\end{equation}
\end{lemma}

\begin{proof}
We have
\begin{align*}
P_n&:=\P\bigl(Z_+(n)+w_{n-\gamma_+(n),n}\in [x,x+\Delta),\: \tau_1\le n\bigr)\\
&=\sum_{k=1}^\infty\P\bigl(T_k=n,\: Z_k\in [x,x+\Delta)\bigr)\\[-2pt]
&\quad\strut+\sum_{k=1}^\infty\,\sum_{i=1}^n\,\int_0^\infty \P\bigl(T_k=n-i,\:
Z_k+y\in [x,x+\Delta),\:\tau_{k+1}\ge i+1,\:u_{k+1,i}\in dy\bigr).
\end{align*}
Since  
$\P(\tau\ge 1,\: u_0=0)=1$ and since $(\tau_{k+1},u_{k+1,i})$ and  $(T_k,Z_k)$ are independent for each $k$,
\begin{align*}
P_n&=\sum_{k=1}^\infty\P\bigl(T_k=n,\: Z_k\in [x,x+\Delta)\bigr)\P(\tau\ge 1,\:
u_0=0)\\ &\quad\strut+\sum_{k=1}^\infty\,\sum_{i=1}^n\,\int_0^\infty
\P\bigl(T_k=n-i,\: Z_k+y\in [x,x+\Delta)\bigr)\P(\tau\ge i+1,\:u_i\in dy)\\
&=Q\sum_{k=1}^\infty\,\sum_{i=0}^n\,\int_0^\infty \P\bigl(T_k=n-i,\: Z_k+y\in
[x,x+\Delta)\bigr) \frac{1}{Q}\P(\tau\ge i+1,\: u_i\in dy)\\
&=Q\sum_{k=1}^\infty\,\sum_{i=0}^n\,\int_0^\infty \P\bigl(T_k=n-i,\: Z_k+y\in
[x,x+\Delta)\bigr)\P(\tau^*= i,\: \zeta^*\in dy)\\ &= Q\sum_{k=1}^\infty
\P\bigl(T^*_k=n,\: Z^*_k\in [x,x+\Delta)\bigr)\\ &=Q \P\bigl(Z_+^*(n)\in
[x,x+\Delta),\: \gamma_+^*(n)=0,\: \nu_+^*(n)\ge 1\bigr)\\ &=Q \P\bigl(Z_+^*(n)\in
[x,x+\Delta),\: \gamma_+^*(n)=0\bigr)-Q \P\bigl(Z_+^*(n)\in [x,x+\Delta),\:
\gamma_+^*(n)=0,\: \nu_+^*(n)=0\bigr)\\ &=Q \P\bigl(Z_+^*(n)\in [x,x+\Delta),\:
\gamma_+^*(n)=0\bigr),
\end{align*}
where the last equality follows from the fact that, for $x>0$,
$$
\P(Z_+^*(n)\in [x,x+\Delta),\nu_+^*(n)=0)=0.\qed
$$
\end{proof}

In the particular case where $v$ has an arithmetic distribution, we $\Delta=1$ in Lemma~\ref{l10} and take $x\in\mathbb{Z}$ to obtain the following corollary of Lemma \ref{l10}.

\begin{lemma}\label{l11}
Assume that conditions $S_{\mathrm{I}}$--$S_{\mathrm{IV}}$ hold and let $v$ satisfy condition $[\mathbf{Z}]$. Then, for any integers $n\ge 2$ and $x\ge 1$, we get 
\begin{equation}\label{1.6.}
\P\bigl(Z_+(n)+w_{n-\gamma_+(n),n}=x,\: \tau_1\le n)=Q\P\bigl(Z_+^*(n)=x,\:
\gamma_+^*(n)=0\bigr).
\end{equation}
\end{lemma}

We continue now with the proof of Theorem \ref{th2}.

I. Consider the arithmetic case first. In order to apply results from the paper \cite{Mog1}, we need to introduce, in addition to CRP $Z_+(n)$, another CRP $Z(n)$, since the main results in \cite{Mog1} are obtained for CRP $Z(n)$).

For $n\ge 1$, let 
$$
\nu(n):=\max\{k\ge 1:\: T_k<n\},\qquad \gamma(n):=n-\nu(n).
$$
Then 
$$
Z(n):=Z_{\nu(n)}.
$$
It is easy to see
(we assume that the processes  $Z(n)$ and $Z_+(n)$ are constructed on a common probability space based on the same sequence $\{(\tau_k,\zeta_k)\}$) that, for any $n\ge 1$,
\begin{equation}\label{13.11.1}
\nu_+(n)=\nu(n+1),\qquad Z_+(n)=Z(n+1),\qquad \gamma_+(n)=\gamma(n)-1.
\end{equation}
In particular, the defect $\gamma(n)$  takes values $\{1,2,\ldots\}$ and the defect $\gamma_+(n)$ takes values $\{0,1,2,\ldots\}$. We will use similar notation for CRP 
$\{(\tau^*_k,\zeta^*_k)\}$ and the corresponding functionals, with adding an extra upper-case $``*''$, for example:
$$
\nu^*(n),\quad \nu_+^*(n),\quad Z^*(n),\quad Z_+^*(n),\quad \gamma^*(n),\quad
\gamma_+^*(n),\quad \text{etc.}
$$
The following relations in the domain of the normal and moderately large deviations follow from formulae  
\eqref{13.11.1} of Theorem 2.1, Corollary 2.1 and Theorem 2.1$^*$ of the paper \cite{Mog1} in the case where $x\in\mathbb{N}$, $x-na=o(n)$: as $n\to\infty$, 
\begin{align}
\P(Z(n)=x)\sim \P(Z_+(n)=x) &\sim \P(Z^*(n)=x)\sim \P(Z_+^*(n)=x)\sim
\frac{1}{\sigma\sqrt{2\pi n}}e^{-nD(\frac{x}{n})},\label{1.10}\\
\P\bigl(Z_+^*(n)=x,\: \gamma_+^*(n)=0\bigr)&\sim \P\bigl(Z_+(n)=x,\:
\gamma_+(n)=0\bigr)\sim \frac{1}{{\bf E}\tau}\P(Z_+^*(n)=x).\label{1.11}
\end{align}
Here equivalences \eqref{1.10} show that, under our assumptions, differences between processes $Z(n)$ and $Z_+(n)$, as well as
diffecences related to inhomogeneity disappear for local theorems in the domain of normal and moderately large deviations.

Applying Lemma \ref{l11} and noticing that $Q={\bf E}\tau$ and that, for some $h>0$,
\begin{align*}&
\P\bigl(Z_+(n)+w_{n-\gamma_+(n),n}=x\bigr)\\ &=
\P\bigl(Z_+(n)+w_{n-\gamma_+(n),n}=x,\: \tau_1>n\bigr)+
\P\bigl(Z_+(n)+w_{n-\gamma_+(n),n}=x,\: \tau_1\le n\bigr) \\
&=\P\bigl(Z_+(n)+w_{n-\gamma_+(n),n}=x,\: \tau_1\le n\bigr) +O(e^{-nh}),
\end{align*}
we obtain from 
\eqref{1.10} and \eqref{1.11} the statement of part I of Theorem \ref{th2}:
$$
\P\bigl(Z_+(n)+w_{n-\gamma_+(n),n}=x\bigr)\sim \frac{1}{\sigma\sqrt{2\pi
n}}e^{-nD(\frac{x}{n})}.
$$

II. Our proof of part II of Theorem \ref{th2} is based on the integro-local theorem in the domains of the normal, moredately large and large deviations obtained by Stone~\cite{St}. We formulate this theorem in notation that is convenient to us. We define the rate function for random vector $(\tau,\zeta)$ as 
$$
\Lambda(\theta,\alpha):=\sup_{\lambda,\mu}\{\lambda\theta+\mu \alpha -
A(\lambda,\mu)\}.
$$
Next, we denote by $|\Lambda''(\theta,\alpha)|$ the determinant of matrix 
$\Lambda''(\theta,\alpha)$ of the second derivatives of the rate function $\Lambda(\theta,\alpha)$.

{\renewcommand{\thetheorem}{\arabic{theorem}\, \rm\cite{St}}
\begin{theorem}%\label{th3}
Let the distribution of vector $(\tau_1,\zeta_1)$ coincide with the distribution of $(\tau,\zeta)$ and let conditions $S_{\mathrm{III}}$ and $[{\bf ZR}]$ hold. Then, for some
$\delta>0$ and for a certain sequence $\Delta^{(0)}:=\Delta^{(0)}_n>0$ that tends to zero as $n\to\infty$, and for any 
$(x,y)\in\mathbb{Z}\times\mathbb{R}$ such that $|\theta-a_\tau|+
|\alpha-a_\zeta|\le\delta$, where
$(\theta,\alpha):=\Bigl(\nfrac{x}{n},\nfrac{y}{n}\Bigr)$, the following equality holds:
$$
\P\bigl(T_n=x,\: Z_n\in [y,y+\Delta)\bigr)= \frac{\Delta
\sqrt{|\Lambda''(\theta,\alpha)|}}{2\pi n} \exp\{-n\Lambda(\theta,\alpha)\}(1+o(1)),
$$
where $\Delta:=\Delta_n\ge\Delta^{(0)}_n$ and $\Delta_n\to 0$ as $n\to\infty$, and the remaining term 
 $o(1)=\varepsilon_n(x,y)$ satisfies the relation 
$$
\lim_{n\to\infty}\,\sup_{\substack{(x,y)\in\mathbb{Z}\times\mathbb{R}\\
|\theta-a_\tau|+ |\alpha-a_\zeta|\le\delta}}|\varepsilon_n(x,y)|=0.
$$
\end{theorem}}
Note that, by applying statements  
$S_{\mathrm{I}}$--$S_{\mathrm{IV}}$ and repeating all phases of the proofs of Theorems~2.1 and  2.1$^*$ and of Corollary 2.1 from \cite{Mog1}, we obtain natural analogues of these statements where symbols $=x$ are replaced by $\in
[x,x+\Delta)$ and an additional coefficient $\Delta$ appears in the right-hand side. The rest of the proof of part II completely repeats the corresponding piece of the proof of part I.\qed

\smallskip
We formulate and prove now the last results. 

\begin{theorem}\label{th4}
The characteristics  
$a$, $\sigma$ and $D(\alpha)$ that are used in Theorem 1 do not depend on a choice of admissible constants $c_1$ and $c_2$.
\end{theorem}

\begin{proof}
We have shown that the local theorems in the domains of the normal and moderately large deviations for processes $w_{0,n}$ and $Z_+(n)$ look identical, and the formulation of Theorem 1 includes characteristics $a$, $\sigma$ and $D(\alpha)$ that are uniquely determined by vector $(\tau,\zeta)$ (see \eqref{17.11.1}--\eqref{17.11.3}) that ``drives'' CRP $Z_+(n)$. Since  random vector $(\tau,\zeta)$ is defined in terms of arbitrarily chosen constants $c_1$ and $c_2$ that satisfy condition 
\eqref{defc}, one can guess that the characteristics $a$, $\sigma$ and  $D(\alpha)$ may depend on these constants too. However, we will show now that this is  not the case. Namely, we will show that,
for any other pair of constants $\widetilde{c}_2\le\widetilde{c}_1$ satisfying 
\eqref{defc}, the corresponding characteristics $\widetilde{a}$, $\widetilde{\sigma}$ and $\widetilde{D}(\alpha)$ determined by vector
$(\widetilde{\tau},\widetilde{\zeta})$ coincide, in fact, with $a$, $\sigma$ and~$D(\alpha)$.

Based on what have we proved already, one can see that, for the processes $w_{0,n}$ and $\widetilde{Z}_+(n)$, a similar local theorem holds, with corresponding characteristics $\widetilde{a}$, $\widetilde{\sigma}$
and~$\widetilde{D}(\alpha)$. This means that, for the process $w_{0,n}$, two local theorems take place, and their statements differ in characteristics $a,\sigma,D(\alpha)$ and 
$\widetilde{a},\widetilde{\sigma},\widetilde{D}(\alpha)$ only. The local theorems in the domain of the normal deviations clearly lead to the corresponding laws of large numbers: for any $\varepsilon>0$, one has 
$$
\lim_{n\to\infty}\P(|w_{0,n}-an|\le n\varepsilon)=1,\qquad
\lim_{n\to\infty}\P(|w_{0,n}-\widetilde{a}n|\le n\varepsilon)=1.
$$
Then $\widetilde{a}=a$, with necessity.

Next, we may rewrite the statements of these local theorems using equality $\widetilde{a}=a$ and arrive at the equivalence 
$$
\frac{1}{\sigma\sqrt{n}}e^{-nD(\frac{x}{n})}\sim
\frac{1}{\widetilde{\sigma}\sqrt{n}}e^{-n\widetilde{D}(\frac{x}{n})},
$$
that is valid for any
$x=x_n\in\mathbb{Z}$ in the area $\Bigl|\nfrac{x}{n}-a\Bigr|=o(1)$. Again, necessarily, the latter equivalence implies equality $\sigma=\widetilde{\sigma}$ and identity of analytic (in a neighbourhood of point 
$\alpha=a$) functions $\widetilde{D}(\alpha)=D(\alpha)$.\qed
\end{proof}

\section{Auxiliary result}

In our proofs, we use an auxiliary result that is not original. We include its proof because it is very short.

\begin{lemma}\label{l4}
Let $S_n=\sum\limits_{i=1}^n X_i$\textup, $n=1,2,\ldots\strut$ be a sequence of independent and identically distributed non-negative random variables 
$\{X_i\}$ and let $N$  be a counting random variable.  Assume that
$$
\E e^{C X_1}<\infty \ \text{and} \ \E e^{C N}<\infty. \quad \text{for some}\ C>0.
$$
Then one can choose a constant 
$b>0$ such that 
$$
\E e^{b S_N}<\infty.
$$
\end{lemma}

\begin{proof}
Take any $a> \E X_1$. Then 
\begin{equation}\label{08.01.1}
S_N=\sum_{i=1}^N (X_i-a)+a N\le\sup_{n\ge 0} (S_n-na)+aN \equiv R+aN,
\end{equation}
here we let $S_0=0$. From \eqref{08.01.1} and from the elementary inequality $e^{x+y}\le e^{2x}+e^{2y}$, for any $b>0$ we have 
$$
e^{bS_N}\le e^{2bR}+ e^{2baN}\le 1+ \sum_{n=1}^{\infty} e^{2b(S_n-na)}+e^{2baN}
$$
and, therefore, 
\begin{equation}\label{randomsum}
\E e^{bS_N}\le 1+\sum_{n=1}^{\infty}\bigl(\E e^{2b(X_1-a)}\bigr)^n+ \E e^{2baN}.
\end{equation}
Since $a>\E X_1$, then one can choose $b>0$ so small that $2b\max (1,a)<C$ and $\E e^{2b(X_1-a)}<1$. For such $b$, the right-hand side of \eqref{randomsum} is finite too.\qed
\end{proof}

\section*{FUNDING}

The paper is supported by the Mathematical Center in Akademgorodok, grant 075-15-2019-1675 by the Ministry of Science and Higher Education (A.V. Logachov, A.A. Mogulskii and S.G.Foss); and also by a joint Russian-French grant RFBR 19-51-15001 -- CNRS 193-382 (S.G. Foss and T. Konstantopoulos).

%% file: _2.bbl
\begin{thebibliography}{12}

\bibitem{Coh}
Cohen, J.E., Briand, F., and Newman, C.M., \emph{Community Food Webs: Data and
Theory}, Berlin: Springer, 1990.

\bibitem{New1}
Newman, C.M., Chain Lengths in Certain Random Directed Graphs, \emph{Random
Structures Algorithms}, 1992, vol.~3, no.~3, pp.~243--253.
\url{https://doi.org/10.1002/rsa.3240030304}

\bibitem{Gel}
Gelenbe, E., Nelson, R., Philips, T., and Tantawi, A., An Approximation of the
Processing Time for a Random Graph Model of Parallel Computation, in \emph{Proc.\
1986 ACM Fall Joint Computer Conf.\ (ACM'86)}, Los Alamitos, CA: IEEE Computer
Society Press, 1986, pp.~691--697.
\url{https://dl.acm.org/doi/proceedings/10.5555/324493}

\bibitem{New2}
Isopi, M. and Newman, C.M., Speed of Parallel Processing for Random Task Graphs,
\emph{Comm.\ Pure Appl.\ Math.}, 1994, vol.~47, no.~3, pp.~361--376.
\url{https://doi.org/10.1002/cpa.3160470307}

\bibitem{Foss3}
Foss, S. and Konstantopoulos, T., Extended Renovation Theory and Limit Theorems for
Stochastic Ordered Graphs, \emph{Markov Process.\ Related Fields}, 2003, vol.~9,
no.~3, pp.~413--468.

\bibitem{Foss2}
Denisov, D., Foss, S., and Konstantopoulos, T., Limit Theorems for a Random Directed
Slab Graph, \emph{Ann.\ Appl.\ Probab.}, 2012, vol.~22, no.~2, pp.~702--733.
\url{https://doi.org/10.1214/11-AAP783}

\bibitem{Mog1}
Mogul'ski\u\i, A.A. and Prokopenko, E.I., Local Theorems for Arithmetic
Multidimensional Compound Renewal Processes under Cram\'er's Condition, \emph{Mat.\
Tr.}, 2019, vol.~22, no.~2, pp.~106--133 [\emph{Siberian Adv.\ Math.}\ (Engl.\
Transl.), 2020, vol.~30, no.~4, pp.~284--302].
\url{https://doi.org/10.1134/S1055134420040033}

\bibitem{Mog2}
Mogulskii, A.A. and Prokopenko, E.I., The Rate Function and the Fundamental Function
for Multidimensional Compound Renewal Process, \emph{Sib.\ Elektron.\ Mat.\ Izv.},
2019, vol.~16, pp.~1449--1463. \url{https://doi.org/10.33048/semi.2019.19.100}

\bibitem{Foss1}
Foss, S., Martin, J.B., and Schmidt, P., Long-Range Last-Passage Percolation on the
Line, \emph{Ann.\ Appl.\ Probab.}, 2014, vol.~24, no.~1, pp.~198--234.
\url{https://doi.org/10.1214/13-AAP920}

\bibitem{Foss4}
Foss, S. and Konstantopoulos, T., Limiting Properties of Random Graph Models with
Vertex and Edge Weights, \emph{J.~Stat.\ Phys.}, 2018, vol.~173, no.~3--4,
pp.~626--643. \url{https://doi.org/10.1007/s10955-018-2080-3}

\bibitem{Tes}
Tesemnikov, P.I., On the Asymptotics for the Minimal Distance Between Extreme
Vertices in a Generalised Barak--Erd\H{o}s Graph, \emph{Sib.\ Elektron.\ Mat.\ Izv.},
2018, vol.~15, pp.~1556--1565.

\bibitem{St}
Stone, C., On Local and Ratio Limit Theorems, \emph{Proc.\ 5th Berkeley Symp.\ on
Mathematical Statitics and Probability, Univ.\ of California, Berkeley, 1965--66}, Le
Cam, L.M. and Neyman, J., Eds., Berkeley, CA: Univ.\ of California Press, 1967,
vol.~2: Contributions to Probability Theory, Part~2, pp.~217--224.
\url{https://doi.org/10.1525/9780520325340-017}

\end{thebibliography}
